\documentclass{jfm}

\usepackage{graphicx}
\usepackage{newtxtext}
\usepackage{newtxmath}
\usepackage{natbib}
\usepackage{hyperref}
\hypersetup{
    colorlinks = true,
    urlcolor   = blue,
    citecolor  = black,
}

\newcommand{\RomanNumeralCaps}[1]

\usepackage[]{algorithm2e}
\usepackage{booktabs}

\newcommand{\bn}{\boldsymbol{n}}
\newcommand{\bu}{\boldsymbol{u}}
\newcommand{\bv}{\boldsymbol{v}}
\newcommand{\bx}{\boldsymbol{x}}
\newcommand{\bW}{\boldsymbol{W}}
\newcommand{\bsigma}{\boldsymbol{\sigma}}
\newcommand{\blambda}{\boldsymbol{\lambda}}
\newcommand{\dd}{\mathrm{d}}
\newcommand{\RR}{\mathbb{R}}

\newcommand{\half}{\mbox{$\frac{1}{2}$}}

\DeclareMathOperator*{\Span}{span}

\title{Numerical approximation of viscous contact problems applied to glacial sliding}

\author{Gonzalo G. de Diego \corresp{\email{gonzalezdedi@maths.ox.ac.uk}},
  Patrick E. Farrell,
 \and Ian J. Hewitt}

\affiliation{Mathematical Institute, University of Oxford, Oxford, OX2 6GG, UK}

\begin{document}
\maketitle

\begin{abstract}
	Viscous contact problems describe the time evolution of fluid flows in contact with a surface from which they can detach and reattach. These problems are of particular importance in glaciology, where they arise in the study of grounding lines and subglacial cavities. In this work, we propose a novel numerical method for solving viscous contact problems based on a mixed formulation with Lagrange multipliers of a variational inequality involving the Stokes equation. The advection equation for evolving the geometry of the domain occupied by the fluid is then solved via a specially-built upwinding scheme, leading to a robust and accurate algorithm for viscous contact problems. We first verify the method by comparing the numerical results to analytical results obtained by a linearised method. Then, we use this numerical scheme to reconstruct friction laws for glacial sliding with cavitation. Finally, we compute the evolution of cavities from a steady state under oscillating water pressures. The results depend strongly on the location of the initial steady state along the friction law. In particular, we find that if the steady state is located on the downsloping or rate-weakening part of the friction law, the cavity evolves towards the upsloping section, indicating that the downsloping part is unstable.
\end{abstract}

\begin{keywords}
\end{keywords}

\section{Introduction}

Viscous contact problems are time-dependent fluid flow problems in which the fluid is in contact with a solid surface from which it can detach and reattach. Contact problems of this type arise when modelling glacial ice flow, which is typically treated as a viscous fluid flow \citep{schoof2013}. On a large scale, they are relevant to marine ice sheets with a grounding line \citep{schoof2007, schoof2011} and, on a smaller scale, to the formation of subglacial cavities when the ice slides over bedrock undulations \citep{fowler1986, schoof2005, gagliardini2007}. These problems share a very similar mathematical structure and are of great importance for understanding ice sheet dynamics and predicting future sea level rise. 

In this paper, we present a novel numerical approach for solving viscous contact problems in an accurate and robust way. This method relies on the formulation of the contact problem as a variational inequality. There exist many different approaches to solving variational inequalities with finite element methods. Of these, we find that solving the contact problem with a piecewise constant Lagrange multiplier is particularly well suited since it allows us to satisfy a discrete version of the contact boundary conditions exactly in the discretisation. We find that this latter property of the numerical scheme enables us to evolve the surface of the viscous flow in a robust way. To the best of our knowledge, viscous contact problems have only previously been solved as variational inequalities in \citet{stubblefield2021}, where a penalty method is used.

Although the numerical method presented in this paper is suitable for any viscous contact problem, here we focus on subglacial cavitation and its application to glacial sliding. Glacial sliding represents a fundamental component of glacier dynamics and is one of the major uncertainties in ice sheet modelling \citep{ritz2015}. It involves a variety of physical mechanisms among which is cavitation. Sliding with cavitation has been studied theoretically by means of linearised approaches \citep{fowler1986,schoof2005} and these investigations are limited to relatively simple scenarios with two-dimensional geometries, smooth bedrocks with small variations and steady conditions. Numerical methods, on the other hand, allow the extension to more complicated situations in a more straightforward manner~\citep{gagliardini2007,helanow2020}. In this work, we attempt to more fully exploit the mathematical structure of the viscous contact problem, to improve the accuracy and robustness of the numerical methods employed.

In Section 2 of this paper, we introduce our algorithm for solving viscous contact problems, using the setup of the glacial cavitation problem to provide a concrete context. We then present two applications of the algorithm. First, in Section \ref{sec:steady-sliding}, we compute the steady sliding law for flow over a sinusoidal bed for linear and nonlinear rheologies, expanding on the results in \citet{gagliardini2007}. Furthermore, in this section we compare our results with those obtained from the linearised theory to validate the algorithm. Then, in Section \ref{sec:unsteady-sliding}, we explore the effects of unsteady water pressures on glacial sliding by calculating the basal sliding velocities and cavity shapes under oscillating water pressures.

\section{Formulation and approximation of the problem of subglacial cavitation}\label{sec:form-approx}

We focus here on formulating the viscous contact problem of subglacial cavitation and describing a finite element scheme to approximate it. The problem described here is the same, though with some different notation, as the subglacial cavitation problem studied in \citet{gagliardini2007}. A linearised version, assuming small amplitude of the bed bump, is equivalent to the problem studied in \citet{fowler1986} and \citet{schoof2005}. Subglacial cavitation is mathematically very similar to the problem of a marine ice sheet with a grounding line. Therefore, the extension of the finite element scheme presented here to grounding line problems requires minor modifications, such as the inclusion of non-zero tangential stress (a friction law), and a free-surface at the ice-atmosphere interface.

\begin{figure}
	\centering
	\includegraphics[width=0.8\textwidth]{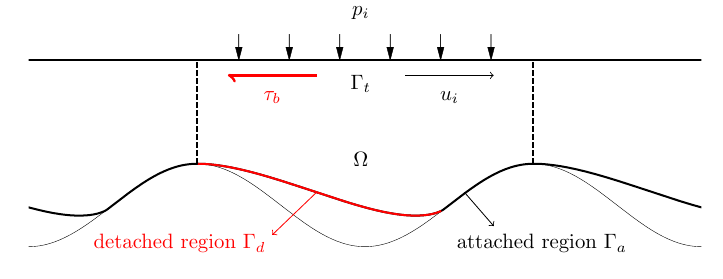}
	\caption{The periodic domain $\Omega$ on which we model the evolution of a subglacial cavity.}
	\label{fig:domain}
\end{figure}

Subglacial cavitation occurs at the ice-bedrock interface, over length scales corresponding to the size of the bedrock obstacles. These length scales are generally several orders of magnitude smaller than those of the glacier. For this reason, the computational domain $\Omega$ in which we model the formation of cavities is a thin layer of ice of finite height located under a larger mass of ice, see Figure \ref{fig:domain}. We assume the bedrock, and therefore also $\Omega$, to be periodic in the horizontal direction. The upper boundary $\Gamma_t$ represents a fictional boundary separating $\Omega$ from the ice above. The height of the bedrock is given by the function $b(x)$ and the height of the lower boundary by $\theta(x,t)$. The lower boundary is divided into an attached region $\Gamma_a$ where $\theta(x, t) = b(x)$, and a detached (i.e.~cavitated) region $\Gamma_d$ where $\theta(x, t) > b(x)$. 

Glacial ice is generally modelled with the Stokes equations. If we denote the velocity and pressure by $\bu = (u,v)$ and $p$ respectively, these can be written as
\begin{subequations}\label{eq:stokes-eq}
	\begin{align}
		- \nabla\cdot\left(2\eta\varepsilon(\bu) \right) + \nabla p &= 0, \label{eq:stokes-eq-1}\\
		\nabla\cdot \bu &= 0.
	\end{align}
\end{subequations} 

Here, $\eta = \eta(|\varepsilon(\bu)|)$ is the effective viscosity of the ice, which is usually modelled with Glen's law \citep{glen1958},
\begin{equation}\label{eq:glens-law}
	\eta(|\varepsilon(\bu)|) = \frac{1}{2}A^{-1/n} \left( \frac{1}{2} |\varepsilon(\bu)|^2 \right)^\frac{1-n}{2n},
\end{equation}
where $n\approx 3$ for ice and $A$ is a potentially temperature-dependent parameter, but which we treat as constant. The tensor $\varepsilon(\bu)$ represents the strain rate of the ice and is given by the symmetric component of the velocity gradient,
\begin{equation}
	\varepsilon(\bu) = \frac{1}{2}(\nabla\bu + \nabla\bu^\top).
\end{equation}
For a given velocity and pressure field, we define the stress tensor $\sigma = \sigma(\bu,p)$ by 
\begin{equation}
	\sigma = 2\eta\varepsilon(\bu) - pI,
\end{equation}
where $I$ is the identity tensor field. If $\bn$ denotes the outwards-pointing normal vector to the boundary of $\Omega$, the normal and tangential stresses at the boundary are defined by
\begin{equation}
	\sigma_{nn} = (\sigma\bn)\cdot \bn \quad \text{and} \quad \bsigma_{nt} = \sigma\bn - \sigma_{nn}\bn.
\end{equation}
The Stokes equations \eqref{eq:stokes-eq} must be complemented with a set of boundary conditions. Along the cavitated region $\Gamma_d$ of the lower boundary, we assume the ice to be in contact with water at a pressure $p_w$, and we prescribe
\begin{equation}\label{eq:cavity-bc}
	\sigma_{nn} = -p_w \quad \text{and} \quad \bsigma_{nt} = 0 \quad \text{on $\Gamma_d$}.
\end{equation}
For the subglacial cavity problem, we assume this water pressure to be uniform along the length of the bedrock because gravity is unimportant on the spatial scales under consideration.

On the attached region $\Gamma_a$, we assume the ice to be lubricated by a thin layer of water connected to the subglacial drainage system. For this reason, we allow the ice to slide freely and set $\bsigma_{nt} = 0$. The possibility of detachment is realised by enforcing the contact boundary conditions
\begin{equation}\label{eq:contact-bc}
	\bu\cdot\bn \leq 0, \quad \sigma_{nn} \leq -p_w \quad \text{and} \quad (\bu\cdot\bn)(\sigma_{nn}+p_w) = 0 \quad \text{on $\Gamma_a$}.
\end{equation}
Notice that here we are enforcing the ice to either remain attached (i.e.~$\bu\cdot\bn = 0$) whenever the compressive normal stress is larger than the water pressure or to have the possibility of detaching ($\bu\cdot\bn < 0$) if the stress equals the water pressure.

On the top boundary $\Gamma_t$, we enforce the boundary condition $\sigma_{nn} = -p_i$, where $p_i$ is the overburden ice pressure. Finally, we close the system with either the Dirichlet boundary condition 
\begin{equation}\label{eq:bc-Dirichlet}
	u = u_i \quad \text{on $\Gamma_t$},
\end{equation}
or the Neumann boundary condition
\begin{equation}\label{eq:bc-Neumann}
	\bsigma_{nt} = \tau_b\, (1,0)^\top\quad \text{on $\Gamma_t$},
\end{equation}
where $\tau_b$ is the basal shear stress. As explained further on in Section \ref{sec:steady-sliding}, we remark that the horizontal velocity $u_i$ is not the same as the basal sliding speed (though they become approximately equal in the linearised theory). However, the value of $\tau_b$ in \eqref{eq:bc-Neumann} coincides with the basal shear stress along the base of the domain (see \eqref{eq:basal} below) due to the overall balance of forces in the horizontal direction.

The Stokes equation must be coupled to a local advection equation which describes the evolution in time of the cavity roof:
\begin{subequations}\label{eq:advection}
	\begin{align}
		\frac{\partial \theta}{\partial t} + \bu \cdot \left( \frac{\partial \theta}{\partial x}, -1 \right)^\top &= 0, \label{eq:advection-1}\\
		\theta &\geq b\label{eq:advection-2}.
	\end{align}
\end{subequations}

Since $\theta$ is time-dependent, so is the domain $\Omega$ and the attached and detached regions $\Gamma_a$ and $\Gamma_d$ along the lower boundary.

Our goal in the rest of this section is to make use of the fact that the instantaneous Stokes problem \eqref{eq:stokes-eq} together with the contact boundary conditions \eqref{eq:contact-bc} is equivalent to a variational inequality. This mathematical structure can be exploited by a finite element discretisation.

\subsection{A Stokes variational inequality}

The finite element method can be used to compute approximate solutions of the Stokes equation with the contact boundary conditions \eqref{eq:contact-bc}. As a first step, we must find a weak formulation by multiplying \eqref{eq:stokes-eq} by suitable test functions and integrating by parts. Weak formulations of partial differential equations are formulated in terms of weak derivatives and the solutions are generally sought in Sobolev spaces \citep{adams2003}. In order to keep the notation simpler, we assume in this section and in Section \ref{subsec:numerical-method} that we enforce the Neumann boundary condition \eqref{eq:bc-Neumann} instead of \eqref{eq:bc-Dirichlet}. In this case, a suitable space of test functions is the vector-valued Sobolev space $V = \bW^{1,\kappa}(\Omega)$ for the velocity and the Lebesgue space $Q = L^{\kappa'}(\Omega)$ for the pressure. Here, $\kappa = 1 + 1/n$, where $n$ is the coefficient in Glen's law, and $\kappa' = \kappa/(\kappa-1)$ is the H\"older conjugate of $\kappa$.

Due to the contact boundary conditions \eqref{eq:contact-bc}, the associated weak formulation of the Stokes equation is a variational inequality. Given the convex set of admissible velocity fields 
\begin{equation}\label{eq:convex_set}
	K = \left\lbrace \bv\in V : \bv\cdot\bn \leq 0 \text{ on $\Gamma_a$} \right\rbrace,
\end{equation} 
the weak formulation of \eqref{eq:stokes-eq} is: find $(\bu,p)\in K\times Q$ such that
\begin{subequations}\label{eq:vi}
\begin{align}
	a(\bu,\bv - \bu) - b(p, \bv - \bu) &\geq f(\bv - \bu) && \forall \bv\in K, \label{eq:vi-1}\\
	b(q,\bu) &= 0 && \forall q\in Q.
\end{align}
\end{subequations}
Here, the operators $a : V\times V\to \RR$ and $b:Q\times V\to\RR$ are defined by
\begin{equation}
	a(\bu,\bv) = \int_\Omega 2 \eta(|\varepsilon(\bu)|)\, \varepsilon(\bu) : \varepsilon(\bv)\,\dd x \quad \text{and} \quad b(q,\bv) = \int_\Omega q (\nabla\cdot\bv)\,\dd x,
\end{equation}
and the linear functional $f:V\to \RR$ by 
\begin{equation}
	f(\bv) = \int_{\Gamma_t} \tau_b (\bv\cdot(1,0)^\top)\,\dd s - \int_{\Gamma_t} p_i (\bv\cdot\bn)\,\dd s - \int_{\Gamma_a\cup\Gamma_d} p_w (\bv\cdot\bn)\,\dd s.
\end{equation}

As we assume that the water pressure $p_w$ is constant along the bedrock, the incompressibility of the flow velocity can be used to rewrite the right hand side as follows:
\begin{equation}\label{eq:rhs_2}
	f(\bv) = \int_{\Gamma_t} \tau_b (\bv\cdot(1,0)^\top)\,\dd s - \int_{\Gamma_t} N (\bv\cdot\bn)\,\dd s,
\end{equation}
where $N = p_i - p_w$ is the effective pressure. Consequently, the solution only depends on $N$, not on the particular values of $p_i$ and $p_w$.

The variational inequality \eqref{eq:vi} is equivalent to the Stokes equation \eqref{eq:stokes-eq} in the sense that if the velocity $\bu$ is at least twice differentiable and the pressure $p$ differentiable, then integration by parts can be performed in order to arrive at \eqref{eq:stokes-eq} with the boundary conditions introduced in the previous section. Of the three conditions in \eqref{eq:contact-bc}, only the kinematic condition $\bu\cdot\bn \leq 0$ is enforced explicitly in the definition of $K$; the remaining conditions are enforced implicitly in \eqref{eq:vi}.

There exist many different approaches for solving a variational inequality with the finite element method \citep{glowinski1981}. For example, in the penalty method, the constraint $\bu\cdot\bn \leq 0$ in the definition of $K$ is enforced via the addition of a nonlinear penalty term in \eqref{eq:vi-1} and the variational inequality is transformed into a variational equality over $\bW^{1,\kappa}(\Omega)$. This is the approach used in \citet{stubblefield2021}. We find that what works best in combination with the advection equation \eqref{eq:advection} is a mixed method with Lagrange multipliers that enforces an average-wise contact condition. We remark that such a condition could also be enforced via a penalty term, but in this case the contact boundary conditions would not be enforced exactly.

The contact boundary conditions \eqref{eq:contact-bc} can be explicitly enforced in \eqref{eq:vi} with a Lagrange multiplier which essentially represents the sum of the normal stress and the water pressure $\sigma_{nn} + p_w$ along the attached region $\Gamma_a$. In order to introduce this reformulation of \eqref{eq:vi}, we need to define a suitable set of admissible Lagrange multipliers. We define the range of the normal trace operator onto $\Gamma_a$ by
\begin{equation}
	\Sigma = \left\lbrace \bv\cdot\bn |_{\Gamma_a} : \bv\in V \right\rbrace.
\end{equation}
The space of admissible Lagrange multipliers is then the convex cone of elements in the dual space $\Sigma'$ satisfying a weak equivalent of the dynamic contact boundary condition $\sigma_{nn} + p_w \leq 0$  in \eqref{eq:contact-bc}:
\begin{equation}
	\Lambda = \left\lbrace \mu\in\Sigma' : \mu(\xi) \geq 0 \quad \text{$\forall \xi\in\Sigma$ such that $\xi\leq 0$ on $\Gamma_a$} \right\rbrace.
\end{equation}

The mixed formulation of the variational inequality \eqref{eq:vi} may be written as: find $(\bu,p,\lambda)\in V\times Q\times \Sigma'$ such that 
\begin{subequations}\label{eq:mixed}
	\begin{align}
		a(\bu,\bv) - b(p, \bv ) - \lambda(\bv\cdot\bn|_{\Gamma_a}) = f(\bv) & && \forall \bv\in V, \label{eq:mixed-1}\\
		b(q,\bu) = 0 & && \forall q\in Q,\label{eq:mixed-2}\\
		\mu(\bu\cdot\bn|_{\Gamma_a}) \geq 0 \quad \forall\mu\in\Lambda, \quad \lambda\in\Lambda \quad \text{and} \quad \lambda(\bu\cdot\bn|_{\Gamma_a}) = 0 &. \label{eq:mixed-3}
	\end{align}
\end{subequations}
Equation \eqref{eq:mixed-1} can be obtained by multiplying \eqref{eq:stokes-eq-1} by a test function $\bv$, integrating by parts, and setting
\begin{equation}\label{eq:lambda-meaning}
	\lambda = \sigma_{nn} + p_w.
\end{equation}
Formulation \eqref{eq:mixed} offers several advantages for finite element approximation in comparison to \eqref{eq:vi}. First, there is no longer any need to build an approximation of the convex set $K$ because we seek the velocity in the space $V$. Second, the conditions \eqref{eq:mixed-3}, which are a weak analogue of the contact boundary conditions \eqref{eq:contact-bc}, can be rewritten as a nonlinear equation at the discrete level, see \eqref{eq:algebraic-3} below. This means that the contact boundary conditions can be enforced exactly by solving the finite element system analogous to \eqref{eq:mixed}. As we explain shortly, this leads to a robust and stable numerical method for solving the complete subglacial cavitation system. 

Under certain conditions on the geometry of the domain and the linear functional $f$, the variational inequality \eqref{eq:vi} has a unique solution. A proof can be found in \citet{dediego2021}, together with an analysis of the finite element approximation of the mixed system \eqref{eq:mixed}.

\subsection{A numerical method for subglacial cavitation}
\label{subsec:numerical-method}

In this section we present a numerical algorithm for solving the complete subglacial cavitation problem which results from coupling the Stokes equation \eqref{eq:stokes-eq} to the time-dependent advection equation \eqref{eq:advection}. We write the discrete counterparts to these equations and explain how they are coupled by deforming the domain according to a contact criterion. The resulting algorithm is summarised in Algorithm \ref{alg:1}. The solver is implemented in Firedrake \citep{rathgeber2016}, using the version available at \citet{zenodo1}. The code for the viscous contact solver presented in Algorithm \ref{alg:1} is openly available, as described in the data availability statement below.

Subglacial cavitation is a time-dependent problem, so its discretisation requires a partition of a given time interval $[0,T]$ into intervals of length $\Delta t$. At each time step $t_k = k \Delta t$ for $k = 0,1,..., N = T/\Delta t$, we must consider a discrete cavity roof $\theta_h^k$ and domain $\Omega^k$ given by
\begin{equation}
	\Omega^k = \left\lbrace (x,y)\in\RR^2 : 0\leq x \leq L, \quad \theta_h^k(x) \leq y \leq H \right\rbrace.
\end{equation}

Additionally, the lower boundary is the union of the attached and detached regions $\Gamma_a^k$ and $\Gamma_d^k$ respectively. These are also time-dependent and are determined at each time step by the contact criterion \eqref{eq:criterion_contact} below.

\begin{algorithm}[H]
 Set $\theta^0$ \;
 \For{$k = 0,1,2,\dots$}{
  Set $\Gamma^{k}_a$ and $\Gamma^{k}_d$ according to \eqref{eq:criterion_contact}\;
  Find $(\bu^k_h, p^k_h, \lambda^k_h)$ such that \eqref{eq:algebraic} holds for $\theta = \theta^k$ and $\Gamma_a = \Gamma^k_a$\;
  Calculate $\theta^{k+1}$ with \eqref{eq:discrete-advection} \;
  If  $\theta_i^{k+1} < b(x_i)$, set $\theta_i^{k+1} = b(x_i)$\;
  Deform mesh\; 
 }
 \caption{Solution procedure for the subglacial cavity problem}
 \label{alg:1}
\end{algorithm}

\subsubsection{Finite element approximation of the Stokes variational inequality}

On each domain $\Omega^k$, we seek a finite element approximation to \eqref{eq:mixed}. To this end, we consider a triangulation of $\Omega^k$ and the finite element spaces $V_h$, $Q_h$ and $\Sigma_h$ given by piecewise continuous quadratic vector fields, piecewise constant scalar fields and piecewise constant scalar functions on $\Gamma_a^k$ respectively. This velocity-pressure pair is known to be inf-sup stable in two dimensions \citep{boffi2013}. Furthermore, as we explain in Section \ref{subsubsec:advection}, working with a piecewise constant Lagrange multiplier enables the use of a very simple upwinding scheme for solving the advection equation \eqref{eq:advection} in a manner consistent with the discrete contact boundary conditions.

We write $V_h = \Span{\{\bv_i\}}_{i=1}^{N_v}$, $Q_h = \Span{\{q_j\}}_{j=1}^{N_q}$ and $\Sigma_h = \Span{\{\mu_k\}}_{k=1}^{N_\mu}$, where $N_v = \dim{V_h}$, $N_q = \dim{Q_h}$ and $N_\mu = \dim{\Sigma_h}$. For the functions $(\bu_h, p_h, \lambda_h) \in V_h\times Q_h\times \Sigma_h$, we denote the respective degrees of freedom (DoFs) in $\RR^{N_v}$, $\RR^{N_q}$ and $\RR^{N_\mu}$ by $\mathbf{u}$, $\mathbf{p}$ and $\blambda$. In order to write an algebraic counterpart of \eqref{eq:mixed-3}, we need the operator 
\begin{equation}
	\boldsymbol{\gamma_n} : \RR^{N_v} \to \RR^{N_\mu}
\end{equation}
that returns the average normal components of a vector $\bv_h\in V_h$ along the edges on $\Gamma_a^k$. The discrete counterpart to \eqref{eq:mixed} is
\begin{subequations}\label{eq:algebraic}
	\begin{align}
		\mathbf{A}(\mathbf{u}) - \mathbf{B} \mathbf{p} - \mathbf{D}\boldsymbol{\lambda} &= \mathbf{f},\label{eq:algebraic-1}\\
		-\mathbf{B}^\top\mathbf{u} &= 0, \label{eq:algebraic-2}\\
		\blambda + \mathbf{C}(\blambda,\mathbf{u}) &= 0 \label{eq:algebraic-3}.
	\end{align}
\end{subequations}
Here, we have introduced the matrices $\mathbf{B}\in\RR^{N_v\times N_q}$ and $\mathbf{D}\in\RR^{N_v\times N_\mu}$, the vector $\mathbf{f}\in\RR^{N_v}$ and the nonlinear operators $\mathbf{A}:\RR^{N_v}\to\RR^{N_v}$ and $\mathbf{C}:\RR^{N_\mu}\times\RR^{N_v}\to \RR^{N_\mu}$. The matrices are given by the elements $\mathbf{B}_{ij} = b(\bv_i,q_j)$ and $\mathbf{D}_{ij} = \int_{\Gamma^k_a} \mu_j (\bv_i \cdot \bn)\,\dd s$ and the vector $\mathbf{f}$ by
\begin{equation}
	\mathbf{f}_i = - \int_{\Gamma_t} p_i(\bv_i\cdot\bn)\dd s - \int_{\Gamma_a\cup\Gamma_d} p_w(\bv_i\cdot\bn)\dd s + \int_{\Gamma_t} \tau_b (\bv_i\cdot(1,0)^\top)\,\dd s.
\end{equation}
The nonlinear operator $\mathbf{A}$ is given by
\begin{equation}
	[\mathbf{A}(\mathbf{u})]_i = \int_\Omega 2 \eta_\epsilon(|\varepsilon(\bu_h)|)\, \varepsilon(\bu_h) :\varepsilon(\bv_i)\,\dd x,
\end{equation}
where $\eta_\epsilon$ is a regularised form of \eqref{eq:glens-law} that prevents infinite viscosity at zero strain rates; it is defined by
\begin{equation}\label{eq:glens-law-reg}
		\eta_\epsilon(|\varepsilon(\bu)|) = \frac{1}{2}A^{-1/n} \left( \frac{1}{2} |\varepsilon(\bu)|^2 + \epsilon^2 \right)^\frac{1-n}{2n},
\end{equation}
for a regularisation parameter $\epsilon > 0$ which we set to $\epsilon = 10^{-2}$. The nonlinear function $\mathbf{C}$ is
\begin{equation}\label{eq:complementarity}
	\mathbf{C}(\blambda,\mathbf{u}) =  \max{\left\lbrace 0, -\blambda+ c (\boldsymbol{\gamma_n}\mathbf{u}) \right\rbrace}
\end{equation}
for an arbitrary $c>0$ which we set to $c = 1$. We find the computations carried out in this paper to be insensitive to the value of $c > 0$ when $c\in [10^{-6},10^6]$. In \eqref{eq:complementarity}, the $\max$ operation is understood to be carried out for each of the elements in the vector $-\blambda + c (\boldsymbol{\gamma_n}\mathbf{u})\in\RR^{N_\mu}$. The use of \eqref{eq:complementarity} is a common way of expressing contact conditions; a particular advantage is that the nonlinear system \eqref{eq:algebraic} can be solved with a semi-smooth Newton method \citep{hintermuller2002}. Moreover, \eqref{eq:algebraic-3} is equivalent to enforcing 
\begin{equation}\label{eq:discrete-contact}
	\boldsymbol{\gamma_n} \mathbf{u} \leq 0, \quad \boldsymbol{\lambda} \leq 0 \quad \text{and} \quad (\boldsymbol{\gamma_n} \mathbf{u})\cdot \boldsymbol{\lambda} = 0,
\end{equation}
a discrete equivalent of \eqref{eq:contact-bc}. By solving \eqref{eq:discrete-contact}, we see that on each edge on $\Gamma_a^k$ we either have that the edge wishes to detach ($\boldsymbol{\gamma_n} \mathbf{u} < 0$) or to remain attached to the bed ($\boldsymbol{\lambda} \leq 0$ and $\boldsymbol{\gamma_n} \mathbf{u} = 0$). At the continuous level, one cannot write \eqref{eq:mixed-3} as in \eqref{eq:algebraic-3} because in general the Lagrange multiplier $\lambda$ is only an element of the dual space $\Sigma'$ and may not even be a function on $\Gamma_a$ \citep{stadler2007}.

\subsubsection{The discrete advection equation and a contact criterion}\label{subsubsec:advection}

\begin{figure}
\centering
\includegraphics[width=0.7\textwidth, trim = 0cm 0cm 0cm 2cm]{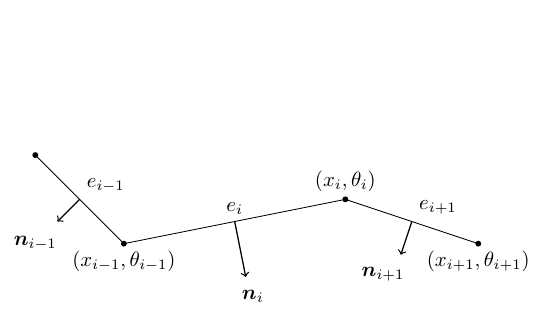}
\caption{Notation and ordering for nodes $\{(x_i,\theta_i)\}$ and edges $\{e_i\}$ along the attached region $\Gamma_a^k$.}
\label{fig:scheme}
\end{figure}

In order to write a discretisation of \eqref{eq:advection}, we introduce some further notation. We denote the points along the lower boundary by $(x_i,\theta_i)$ with the index $i$ increasing from left to right as in Figure \ref{fig:scheme}. The edge between $(x_{i-1},\theta_{i-1})$ and $(x_i,\theta_i)$ is denoted by $e_i$. We write $u^k_{n,i}$ for the value of $\boldsymbol{\gamma_n}\mathbf{u}$ associated to $e_i$ at time $t_k$ (by the definition of $\boldsymbol{\gamma_n}$, this value corresponds to the average value of $\bu_h\cdot\bn$ along $e_i$). 

A numerical scheme that solves the advection equation \eqref{eq:advection-1} should satisfy two conditions. First, it must be stable, as described in standard textbooks \citep{leveque2007} (in Appendix \ref{app:2} we present results with and without stabilisation of \eqref{eq:advection-1}). Additionally, the scheme should be compatible with the discrete contact boundary conditions \eqref{eq:discrete-contact}, in the sense that if for two edges $e_i$ and $e_{i+1}$ we have that $u_{n,i}^k = 0$ and $u_{n,i+1}^k = 0$, then we should have that $\theta^{k+1}_{i} = \theta^{k}_{i}$; that is, the node in between the edges remains unchanged. This means that the discrete counterpart of $\bu\cdot (\partial\theta/\partial x, -1)^\top$ should be defined in terms of $\boldsymbol{\gamma_n} \mathbf{u}$. If not, we generally find that all edges will detach within a few timesteps due to approximation errors. Here, we propose the following explicit scheme:
	\begin{equation}\label{eq:discrete-advection}
		\frac{\theta_i^{k+1} - \theta_i^k}{\Delta t} = \left( \left(\dfrac{\theta^k_i - \theta^k_{i-1}}{x_i - x_{i-1}}\right)^2 + 1 \right)^{1/2} u^k_{n,i},
	\end{equation}
where the velocity at each node is taken from the edge located immediately upstream. This results in an upwinding scheme with the property that penetration cannot occur along $\Gamma_a$, because the average values of $\bu_h\cdot\bn$ along the edges are non-positive, see \eqref{eq:discrete-contact}. However, a node on a previously detached edge can penetrate the bedrock. To avoid this and enforce \eqref{eq:advection-2}, we simply set $\theta^{k+1}_i = b(x_i)$ whenever $\theta^{k+1}_i < b(x_i)$.

Given a cavity roof $\theta_h^k$ at $t_k$, we must set the attached and detached regions $\Gamma_a^{k}$ and $\Gamma_d^{k}$ to solve the Stokes variational inequality at $t_k$. This is carried out by assigning each edge $e_i$ with endpoints $(x_{i-1},\theta^k_{i-1})$ and $(x_i,\theta^k_i)$ to either $\Gamma_a^k$ or $\Gamma_d^k$. We do this by looking at the node downstream of the edge, in accordance with the scheme for the advection equation:
\begin{equation}
\label{eq:criterion_contact}
	\text{$e_i$ belongs to $\Gamma_d^k$ if and only if $\theta^k_i - b(x_i)> \mathrm{tol}$,}
\end{equation}
where the tolerance $\mathrm{tol}$ is set to $10^{-9}$. Consequently, given a cavity roof $\theta^k_h$ at time $t_k$, the attached and detached regions $\Gamma^k_a$ and $\Gamma^k_d$ are fully determined.

Once we have computed a new cavity roof $\theta^{k+1}_h$, the Stokes problem is then solved over a new domain $\Omega^{k+1}$. This means that the mesh is deformed according to $\theta^{k+1}_h$. Since we use vertically extruded meshes for the computations in this paper, a simple algorithm is used which moves each node vertically in a linear manner with respect to the original position of the nodes at $t = 0$. If we denote the mesh nodes over $x_i$ by $y_{ij}$ and let $y = H$ be the upper boundary of the domain, we compute the new vertical positions of the nodes via 
\begin{equation}\label{eq:deform-mesh}
	y_{ij}^{k+1} = \theta^{k+1}_i + \frac{H - \theta^{k+1}_i}{H - \theta^0_i} \left( y^0_{ij} - \theta^0_i \right)
\end{equation}
for a given $y^0_{ij}$, $\theta^0_i$ and $\theta^{k+1}_i$.

\section{Steady sliding with cavitation}\label{sec:steady-sliding}

The sliding of a glacier over its bedrock has been widely studied since Weertman's seminal work in 1957 \citep{weertman1957}. In general, these studies attempt to build a function known as the sliding law that captures the steady relationship between the basal sliding speed $u_b$, the basal shear stress $\tau_b$ and other variables such as the water pressure $p_w$. This sliding law can then be used to prescribe a boundary condition at the ice-bedrock interface in large-scale glacier models which do not resolve the smaller-scale shape of that interface.

As a first application of the algorithm described in the previous section, we build the sliding law for ice flowing over a sinusoidal bed. We first present detailed results for a particular steady state in Section \ref{subsec:steady_cavity} to evaluate the accuracy of the solver and the effect of mesh refinement. Then, in Section \ref{subsec:steady-sliding} we compute sliding laws for different values of the parameter $n$ in Glen's law \eqref{eq:glens-law}.

The computations presented in this section and in Section \ref{sec:unsteady-sliding} have been carried out under the assumption that the water pressure $p_w$ is uniform along the lower boundary. As explained in the previous section, the problem then only depends on the effective pressure $N = p_i - p_w$, rather than individually on $p_i$ or $p_w$ (after subtracting $p_w$ from the stress, the boundary conditions become $\sigma_{nn} = - N$ on $\Gamma_t$, $\sigma_{nn} = 0$ on $\Gamma_d$ and $\sigma_{nn} \leq 0$ on $\Gamma_a$). We non-dimensionalise the calculations by scaling lengths with the wavelength of the bed $L$, velocities with a characteristic velocity scale $U$, and stresses with the characteristic stress scale $(U/(2AL))^{1/n}$. The resulting problem depends only on the non-dimensional amplitude $r$ of the sinusoidal bed, the non-dimensional velocity $u_i$ prescribed along the upper boundary, and the non-dimensional effective pressure $N$.

\subsection{Steady subglacial cavities}\label{subsec:steady_cavity}

\renewcommand{\arraystretch}{1.25}
\begin{table}
\centering
\caption{Information about the different meshes used to compute the steady cavity states together with calculations of the basal shear stress $\tau_b$, the basal sliding speed $u_b$, and the detachment and reattachment points of the cavity.}
\label{tab:steady-state}
\begin{tabular}{cccccc}
\toprule
$n_e$ & mesh cells & $\tau_b$ & $u_b$ & detachment & reattachment \\
\midrule
16 & 96 & 0.014772 & 0.98667 & 0.7500 & 1.0000 \\
32 & 192 & 0.015143 & 0.98633 & 0.7188 & 1.0000 \\
64 & 768 & 0.015484 & 0.98598 & 0.7188 & 1.0000 \\
128 & 3072 & 0.015679 & 0.98577 & 0.7109 & 1.0000 \\
192 & 7296 & 0.015741 & 0.98570 & 0.7135 & 0.9948 \\
\bottomrule
\end{tabular}
\end{table}

\begin{figure}
\centering
\includegraphics[width=\textwidth]{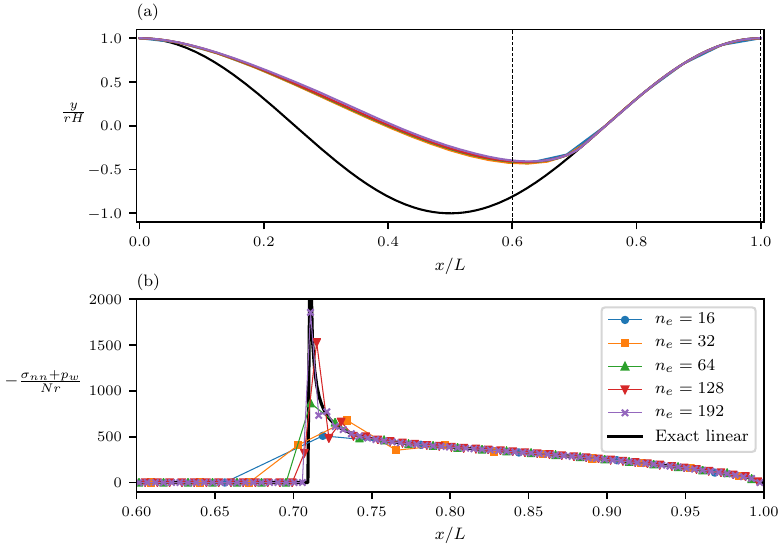}
\caption{(a) Steady cavity shapes and (b) normal stresses along the attached region computed with different mesh sizes. In black, the stress distribution obtained from the linearised theory. Here, the amplitude of the obstacle is set to $r = 0.01$, the scaled effective pressure to $N = 0.3$, and the scaled horizontal velocity at the top to $u_i = 1$.}
\label{fig:steady-state}
\end{figure}

The steady state of equations \eqref{eq:stokes-eq} and \eqref{eq:advection} with the boundary conditions $u = u_i$ and $\sigma_{nn} = -N$ on the top boundary $\Gamma_t$ can be found with Algorithm \ref{alg:1} by evolving the cavity from an initial state until the norm of the discrete derivative $(\theta^{k+1} - \theta^{k})/\Delta t$ is below a prescribed threshold of $10^{-4}$. In this section, we find the steady cavity for a Newtonian flow ($n = 1$ in Glen's law \eqref{eq:glens-law}) over a sinusoidal bed of small amplitude $r = 0.01$ given by $b(x) = rL \cos{\left(2\pi x/L\right)}$. The non-dimensional effective pressure is set to $N = 0.3$ and the non-dimensional horizontal velocity to $u_i = 1$. We first consider a linear rheology and a bedrock of small amplitude in order to compare our results with the analytical solution of the linearised cavitation problem considered in \citet{fowler1986} and \citet{schoof2005}. A brief description of this method is included in Appendix \ref{app:linear}.

We use five different meshes with $n_e$ vertices uniformly distributed along the lower boundary. In Table \ref{tab:steady-state} we present the non-dimensional basal shear stress $\tau_b$ and sliding speed $u_b$ along the cavities. These values are calculated with
\begin{equation}\label{eq:basal}
	\tau_b = - \frac{1}{L} \int_{\Gamma_a\cup\Gamma_d}(\sigma_{nn} + p_w)\,n_x\, \dd s \quad \text{and} \quad u_b = \frac{1}{L}\int_{\Gamma_a\cup\Gamma_d}u\,\dd s,
\end{equation}
where $\bn = (n_x,n_y)$. The expression for $\tau_b$ can be derived from the expression for the force exerted by the ice on the bed as in \cite{schoof2005}. Recall that the Lagrange multiplier $\lambda$ represents the normal stress $\sigma_{nn} + p_w$ along $\Gamma_a$, see \eqref{eq:lambda-meaning}; therefore $\tau_b$ can be calculated from $\lambda$ via 
\begin{equation}\label{eq:basal_lambda}
	\tau_b = - \frac{1}{L} \int_{\Gamma_a}\lambda n_x\, \dd s.
\end{equation}
As we show in the next section, the values of $\tau_b$ and $u_b$ can be used to construct a sliding law.

The formula for $u_b$ presented in \eqref{eq:basal} might seem strange if the subglacial cavity domain is interpreted as a boundary layer between an ice sheet and the bedrock. This would suggest we take the sliding speed to be the average value of the horizontal velocity along the top boundary $\Gamma_t$. However, our computations indicate that, if the height of the domain $H$ is sufficiently large, we can expect the shear stress to approach a constant value and the horizontal velocity to vary with $y^n$ as $y$ approaches $H$, where $n$ here is the exponent in Glen's law \eqref{eq:glens-law}. Therefore, the horizontal velocity along the top boundary depends strongly on the height of the domain. For this reason, and following \citet{gagliardini2007}, we use the equation in \eqref{eq:basal} to calculate $u_b$. In this case, we find that $u_b$ is independent of $H$ for sufficiently large values of $H$. In particular, throughout this paper we set $H=L$. In agreement with \citet{gagliardini2007}, we find this value of $H$ to be sufficiently large.

In Figure \ref{fig:steady-state} we present the steady cavity shape and normal stresses $\sigma_{nn}$ along the attached region $\Gamma_a$ for the different meshes. We can see from these figures that the cavity shape is accurately computed even with the coarsest mesh. Additionally, we also present the stress distribution obtained from the linearised theory, which is uniquely determined for an effective pressure $N$ and a sliding speed $u_b$. The result from the linearised theory plotted in Figure \ref{fig:steady-state} is computed with the value of $u_b$ calculated with the most refined mesh.

The plot for the normal stress distribution demonstrates that the contact conditions \eqref{eq:contact-bc} are satisfied exactly at the discrete level for all of the meshes, because $\sigma_{nn} + p_w \leq 0$. This plot also exhibits the singularity of the normal stresses at the reattachment point. This singularity complicates the approximation of the normal stresses along the attached region and can lead to very inaccurate computations of the sliding law in largely cavitated states. However, Figure \ref{fig:steady-state} also indicates that, with increasing mesh refinement, the solver appears to converge towards the linearised solution. 

At the discrete level, the Lagrange multiplier $\lambda$, which we use to calculate $\tau_b$, is piecewise constant on each edge along the lower boundary. In Figure \ref{fig:steady-state}, these values, the DoFs of $\lambda$, are plotted at the midpoints of each edge. Due to the contact criterion \eqref{eq:criterion_contact}, the edge immediately upstream of the first reattached node is treated as part of $\Gamma_a$ (this is required in order to allow for the possibility of subsequent detachment there), which explains why there are non-zero values of $\sigma_{nn} + p_w$ left of the reattachment node in Figure \ref{fig:steady-state}. This is particularly visible for the coarsest mesh with $n_e = 16$.

\subsection{Computation of the linear and nonlinear steady sliding law}\label{subsec:steady-sliding}

We next perform similar calculations to those of Section \ref{subsec:steady_cavity} but for varying effective pressure $N$ and power-law exponent $n$. This allows us to map out a steady sliding law for ice sliding over a hard bed with cavitation as in \citet{gagliardini2007}. The steady sliding law can be mapped out by only varying $N$ and keeping $u_i$ fixed, since dimensional analysis of the steady problem shows that the scaled basal shear stress $\tau_b/N$ depends only on the ratio $u_b/(ALN^n)$, and not independently on $u_b$ or $N$ \citep{fowler1986} (the same will not be true of the unsteady problem in Section \ref{sec:unsteady-sliding}).

Several previous studies have suggested what form the sliding law should take, both with and without cavitation \citep{kamb1970,fowler1981,fowler1986,gudmundsson1997,schoof2005,gagliardini2007}. The law proposed in \citet{gudmundsson1997} for the uncavitated case can be written as
\begin{equation}\label{eq:nonlinear-sliding-gudmundsson}
	\left(\frac{\tau_b}{rN}\right)^{n} = \alpha(n) \frac{r}{AL} \frac{u_b}{N^n},
\end{equation}
where $r = a/L$ and $\alpha(n)$ is a function depending on $n$. The function $\alpha(n)$ is related to the parameter $c_0$ (which also depends on $n$) considered in \citet{gudmundsson1997} by 
\begin{equation}\label{eq:c0}
	\alpha(n) = \frac{(2\pi)^{n+2}}{2c_0}.
\end{equation}
For a Newtonian flow, the complex analysis method presented in \citet{fowler1986} and \citet{schoof2005} yields an exact solution to the linearised problem. In particular, for high effective pressures, no cavitation occurs and a linear sliding law as in \eqref{eq:nonlinear-sliding-gudmundsson} with $c_0 = 1$ is found. For effective pressures lower than a critical value $8\pi^2 r \eta u_b / L$, cavitation occurs and the sliding law becomes non-linear, varying with $N$ as well as $u_b$.

\begin{figure}
\centering
\includegraphics[width=\textwidth]{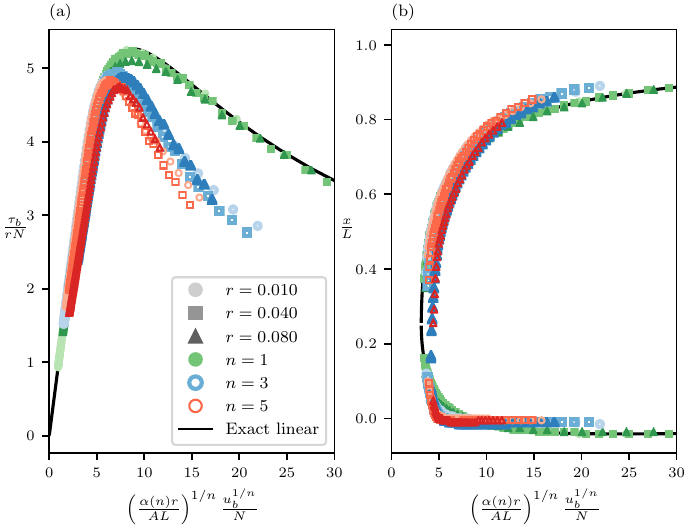}
    \caption{(a) Computed sliding law for steady glacial sliding with cavitation. (b) Cavity endpoints. The parameter $n$ in Glen's law \eqref{eq:glens-law} is set to $n = 1$ (green), 3 (blue) and 5 (red); for each $n$, we compute the sliding law for beds of amplitude $r = 0.01$ (circles, light), 0.04 (squares) and 0.08 (triangles, dark). For these computations, we set the scaled velocity at the top of the domain to $u_i = 1$ and the scaled viscosity parameter to $A = 0.5$. We use a mesh with 192 cells along the lower boundary. The parameter $\alpha(n)$ is computed from the slope of the curve near the origin for the lowest $r$.}
    \label{fig:sliding}	
\end{figure}

We compute the sliding law over a sinusoidal bed of different amplitudes $r$ and for $n=1$, 3 and 5, and we plot the results in Figure \ref{fig:sliding} using the scaling suggested by \eqref{eq:nonlinear-sliding-gudmundsson}. The mesh has 192 vertices along the lower boundary and 7296 cells. The location of the cavity endpoints is also plotted in Figure \ref{fig:sliding}, along with the solution to the linearised problem calculated with the method from \citet{fowler1986} and \citet{schoof2005}. For each $n$, the parameter $\alpha(n)$ is computed by calculating the slope of the curve near the origin (where there is no cavitation) for the lowest value of $r$. The corresponding values of the parameter $c_0$ can then be calculated from \eqref{eq:c0}; these values can be found in Table \ref{tab:c0} together with those obtained in \citet{gudmundsson1997} and \citet{gagliardini2007}. We see that the results obtained in these works are broadly similar to ours and that the value $c_0 = 1$ obtained in the linearised theory is approached in all cases when $n = 1$.

The computed sliding laws with cavitation in Figure \ref{fig:sliding} are multivalued for $\tau_b/(rN)$ as expected \citep{fowler1986,schoof2005}. This aspect of the law justifies the use of the Dirichlet boundary condition $u = u_i$ instead of the alternative Neumann boundary condition. We find that if we use the Neumann boundary condition $\bsigma_{nt} = \tau_b(1,0)^\top$ and initiate the cavity from a fully attached state, the solver always evolves to the steady state associated to the upsloping region of the curve (see also Section \ref{sec:unsteady-sliding}). From Figure \ref{fig:sliding}, we also deduce that the validity of the sliding law \eqref{eq:nonlinear-sliding-gudmundsson} along the linear segment of the curves (where either little or no cavitation has occurred) decreases with increasing values of $r$ and $n$. For example, when $n = 5$, one can observe that the linear segment of the sliding law for $r = 0.08$ clearly does not collapse onto the corresponding linear segment for $r = 0.01$. As soon as the cavity size increases and the sliding laws cease to be linear, the aspect of these curves largely differ for different values of $n$.

In Figure \ref{fig:sliding}, we use a different scaling to the one used in \citet{gagliardini2007}. In \citet{gagliardini2007}, the computed maximum value reached by $\tau_b/N$ is included in the scaling for the sliding law. In this way, the maximum value reached by the scaled sliding law equals 1 by design. However, we preferred the scaling based on \eqref{eq:nonlinear-sliding-gudmundsson} because it contains fewer terms that are unknown a priori. It is also worth mentioning that, for different values of $n$, the curves in Figure \ref{fig:sliding} do not collapse into a single curve when plotted with the scaling from \citet{gagliardini2007}.

For the linear case with $n = 1$, the numerical results computed with the finite element solver highly resemble those obtained with the linearised solution. For $r=0.08$ a slight difference with the linearised solution can be seen near the peak of the sliding law. This difference is probably a consequence of nonlinear effects that are accentuated with increasing amplitudes of bedrock roughness.

\renewcommand{\arraystretch}{1.25}
\begin{table}
\centering
\caption{Value of the parameter $c_0$ associated to the sliding laws. This parameter is computed with the slope of these curves near the origin.}
\label{tab:c0}
\begin{tabular}{c|ccc}
\toprule
 & $r = 0.01$ & \citet{gudmundsson1997} & \citet{gagliardini2007} \\
\midrule
$n = 1$ & 1.0014 & 0.9936 & 0.9771\\
$n = 3$ & 0.3434 & 0.3294 & 0.2769 \\
$n = 5$ & 0.1255 & 0.1153 & - \\
\bottomrule
\end{tabular}
\end{table}

\section{Unsteady sliding}\label{sec:unsteady-sliding}

In the previous section, the sliding law was constructed by computing steady cavity states. However, field measurements from alpine glaciers and from the Greenland Ice Sheet have found short term variations in the water pressure, on timescales down to hours \citep{iken1981, iken1986, sugiyama2004, andrews2014, hoffman2016}. In these studies, variations of water pressure have been correlated with variations in surface speeds, vertical strain and uplift. Subglacial cavitation has been considered a possible mechanism causing these correlations \citep{iken1986, mair2002, sugiyama2004}. These observations motivate an investigation of glacier sliding under unsteady conditions. In this section, we therefore compute the evolution in time of subglacial cavities under oscillating water pressures and calculate the corresponding unsteady basal sliding speed and shear stresses. The study published by Iken of the transient stages between steady cavity shapes \citep{iken1981} is the only numerical investigation of unsteady cavitation solving the Stokes problem known to the authors. 

We initialise the computations from a steady state corresponding to a point in the sliding law determined by an effective pressure $N_0$, a basal sliding speed $u_{b,0}$, and a basal shear stress $\tau_{b,0}$. Instead of prescribing the Dirichlet boundary condition $u = u_i$ on $\Gamma_t$, we enforce the Neumann boundary condition $\bsigma_{nt} = \tau_{b,0} (1,0)^\top$ on $\Gamma_t$. We consider it more physically realistic to have the basal shear stress fixed rather than the sliding speed because we can expect the basal stresses to balance the gravitational driving stresses, which are essentially fixed on these timescales. In practice, if water pressure variations are spatially localised, the driving stress can be transferred to neighbouring regions of the bed, but it is not easy to account for this within the current boundary-layer treatment of the problem. We set $n = 3$ in Glen's law to model the nonlinear rheology of ice. Algorithm \ref{alg:1} is used with a mesh with 192 elements along the lower boundary over a sinusoidal bed of amplitude $r = 0.08$.

\begin{figure}
	\centering
	\includegraphics[width=\textwidth]{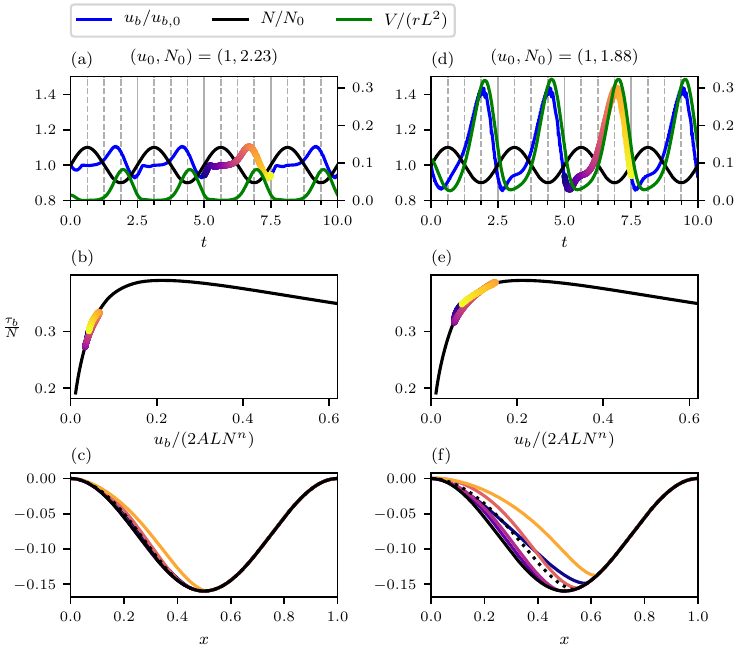}
	\caption{Unsteady cavitation for imposed oscillating effective pressures $N$ and fixed basal shear stress $\tau_b$ around the states with $(u_{b,0}, N_0) = (1, 2.2281)$ (left) and $ (1, 1.8843)$ (right). In (a) and (d), the evolution of the cavity volume $V$ (right axis) and basal sliding speed $u_b$ (left axis, also for $N$) are plotted. One period of each solution is superimposed on the steady sliding law in (b) and (e), as indicated by the coloured dots. In (c) and (f), the cavity shapes are plotted at different time instants with coloured lines; the dotted lines represent the steady cavity shapes for $(u_{b,0},N_0)$.}
	\label{fig:unsteady-1}	
\end{figure}

The effects of unsteady water pressures differ depending on the initial steady state from which we evolve the cavity. To illustrate this, we first evolve two different points along the upsloping component of the steady sliding law by oscillating the effective pressure with an amplitude of $0.1 N_0$ and a fixed non-dimensional frequency of 0.4. As a reference, note that one non-dimensional time unit is approximately the time taken for ice at the top of the domain to traverse one wavelength of the bed. The results are plotted in Figure \ref{fig:unsteady-1}. These results indicate that, with increasing cavitation, the amplitude of the sliding speed increases. For the case of small cavitation (Figure \ref{fig:unsteady-1}(a)), the sliding speed is slightly out of phase with the effective pressure. However, this phase difference disappears with larger cavitation, as observed in Figure \ref{fig:unsteady-1}(d) and also in Figures \ref{fig:unsteady-2} and \ref{fig:unsteady-3} below. This implies that the maximum sliding speed is most often reached when the effective pressure is lowest. Field measurements have also found maximum surface speeds to take place at moments of maximum water pressures \citep{iken1986,sugiyama2004}. On the other hand, the phase difference between the sliding speeds and the cavity volume appears to change in each numerical test: in Figure \ref{fig:unsteady-1}(a), one can observe that the maximum sliding speed is reached when the cavity is still growing, while in Figure \ref{fig:unsteady-1}(d) it is reached at the time of maximum cavitation. There are slight oscillations in the computed sliding speed when the cavity volume is at its largest. These are numerical artefacts due to the stress singularity at the reattachment point of the cavity having an increasing effect on the overall solution of the problem as the cavity volume grows. In these situations, a small displacement of the reattachment point has a large effect on the stress distribution along the bed and therefore also on the computed sliding speed.

\begin{figure}
	\centering
	\includegraphics[width=\textwidth]{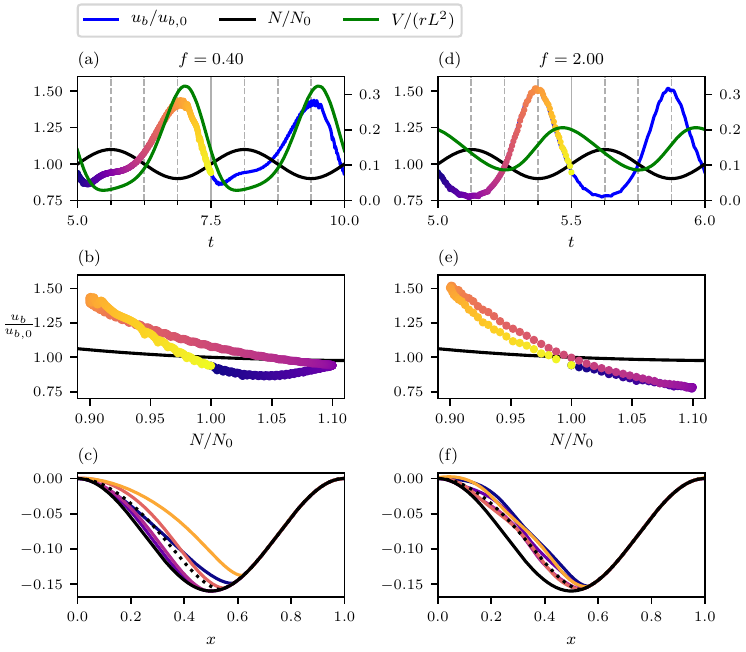}
	\caption{Unsteady cavitation around the state with $(u_{b,0}, N_0) = (1, 1.8843)$ for imposed oscillating effective pressures $N$, with non-dimensional frequencies of $0.4$ (left) and $2$ (right), and fixed basal shear stress $\tau_b$ . In (a) and (d), the evolution of the cavity volume $V$ (right axis) and basal sliding speed $u_b$ (left axis, also for $N$) are plotted. The effective pressure is plotted against the sliding speed in (b) and (e) in coloured dots together with the steady solution is black. In (c) and (f), the cavity shapes are plotted at different time instants with coloured lines; the dotted lines represent the steady cavity shapes for $(u_{b,0},N_0)$.}
	\label{fig:unsteady-2}
\end{figure}

As mentioned above, the non-dimensional time it takes a fluid particle to traverse the domain is of order $t\approx 1$. Therefore, the scaled frequency $f = 0.4$ can be considered a relatively slow frequency that allows the cavity to approximately follow the steady shapes associated to the effective pressure at each instant in time as calculated in Section \ref{sec:steady-sliding}. In Figure \ref{fig:unsteady-2}, we compare results obtained with frequencies $f = 0.4$ and $f = 2$ to examine the effect of faster oscillations in the water pressure. With a higher frequency, the magnitude of the change of cavity volume is significantly lower. In contrast, the amplitude of the sliding speed increases slightly with a higher frequency. The phase difference between the velocity and the cavity volume also changes when the frequency is increased. For a high frequency, the maximum velocity is no longer attained when the cavity reaches its largest extent, but before, when the cavity is still growing. More specifically, Figure \ref{fig:unsteady-2}(d) indicates that the maximum of the sliding speed and of the rate of change of the cavity volume occur simultaneously.

\begin{figure}
	\centering
	\includegraphics[width=\textwidth]{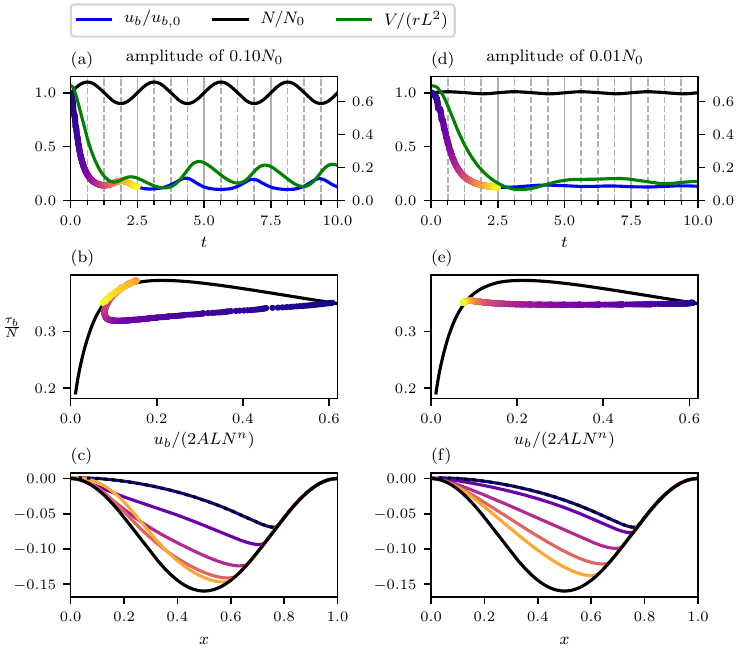}
	\caption{Unsteady cavitation around the state with $(u_{b,0}, N_0) = (1, 1.0937)$ for imposed oscillating effective pressures $N$, with amplitudes of $0.1N_0$ (left) and $0.01N_0$ (right), and fixed basal shear stress $\tau_b$ . In (a) and (d), the evolution of the cavity volume $V$ (right axis) and basal sliding speed $u_b$ (left axis, also for $N$) are plotted. One period of each solution is superimposed on the steady sliding law in (b) and (e), as indicated by the coloured dots. In (c) and (f), the cavity shapes are plotted at different time instants with coloured lines; the dotted lines represent the steady cavity shapes for $(u_{b,0},N_0)$.}
	\label{fig:unsteady-3}
\end{figure}

We can therefore expect a phase difference between the sliding speed and the cavity volume to arise for large frequencies. For very low frequencies, at each time instant $t$, the cavity shape is approximately that of the corresponding steady state for the values of $N(t)$, $u_b(t)$ and $\tau_b(t)$. Therefore, since $u_b$ increases with the cavity volume under steady conditions, this phase difference will disappear. When comparing the panels (a) and (d) of Figure \ref{fig:unsteady-1}, we see that this phase difference is larger when the cavity is smallest. This could indicate that the characteristic time scale with respect to which we measure the frequency increases for smaller cavity volumes.

Following \citet{sugiyama2004}, in Figure \ref{fig:unsteady-2} we also plot the sliding speed against the effective pressure throughout one cycle. For the lower frequency $f = 0.4$, a clear loop arises in which the sliding speed is greater during cavity growth. For the higher frequency $f = 2$, the loop nearly collapses into a single line. These plots can be compared with those obtained from field measurements in \citet[Figure 5]{sugiyama2004}; a qualitative similarity between both is that higher speeds are reached when $N$ decreases. Additionally, similar plots are presented in \citet[Extended Data Figure 4]{andrews2014}. These plots also show the extent to which unsteady sliding can differ from its steady counterpart.

The downsloping section of the friction curve produces a so-called rate-weakening sliding regime in which, for a supposed fixed effective pressure, an increase in the sliding speed is accompanied by a decrease in the basal drag. Rate-weakening sliding has been observed in a laboratory setting for ice sliding over a sinusoidal bed \citep{zoet2015}, although several authors have questioned whether such a sliding regime can arise for more realistic bed geometries \citep{fowler1987,schoof2005, helanow2021}. An implication of rate-weakening sliding is that the sliding law becomes double-valued, as seen in Figure \ref{fig:sliding}. This invalidates the commonly used shallow ice approximation of the Stokes equations which requires the friction law to be invertible \citep{schoof2013}.

In Figure \ref{fig:unsteady-3} we perturb a steady state along the downsloping section of the curve with an oscillating effective pressure of nondimensional frequency $f = 0.4$ with the amplitude set to $0.1N_0$ (left panels) and $0.01N_0$ (right panels). As shown in Figures \ref{fig:unsteady-3}(b) and \ref{fig:unsteady-3}(d), we observe that, for perturbations with both large and small amplitudes, the cavity quickly evolves towards the steady state along the upsloping section for a similar value $\tau_{b,0}/N_0$. In fact, we find that this phenomenon continues to occur for different frequencies in the oscillations of the effective pressure and different steady states along the downsloping section of the sliding law. This finding could offer an additional reason to not use a sliding law with a rate-weakening regime: since such a regime is unstable in the sense described above, we could expect it to be unachievable under natural conditions.

\section{Discussion and conclusions}

In this work, we have presented a novel numerical method for solving viscous contact problems by formulating the Stokes variational inequality as a mixed problem with Lagrange multipliers. We have also introduced a numerical scheme for the advection equation specifically designed to ensure both stability and coherence with the contact conditions of the Stokes problem. This combination leads to a robust and simple method for solving viscous contact problems. Although we have only applied our method to the problem of subglacial cavitation, it can be extended to the ice-ocean grounding line problem. In Section \ref{sec:steady-sliding} we have used this method to compute steady cavity configurations. In particular, we have validated the numerical method by comparing our results over bedrocks with small amplitudes to the linearised approach in \citet{fowler1986} and \citet{schoof2005} and we have reconstructed steady sliding laws for different values of the parameter $n$ in Glen's law \eqref{eq:glens-law}. 

Finally, in Section \ref{sec:unsteady-sliding}, we have explored the temporal evolution of cavities under unsteady effective pressures and its effect on glacial sliding. One of the features of unsteady sliding studied in this work is the phase difference between the sliding speed, the cavity volume and the effective pressures. Our results show that with increasing frequencies, the phase difference between the sliding speed and cavity volume increases. They also seem to indicate that the maximum sliding speed occurs at the point of minimum effective pressure, at least for sufficiently cavitated states. Similar phase differences have been found in data obtained from field measurements \citep{iken1986, sugiyama2004, andrews2014}. Although our results could offer an explanation in terms of an idealised model, it should be noted that changes in measurements of surface elevation of an ice sheet can be the result of many cavities in different states. Finally, we also find that when we fix the value of $\tau_b$ as a Neumann boundary condition, the downsloping section of the sliding law, also known as the rate-weakening regime, is unstable under finite perturbations. In particular, if we perturb a steady state along the downsloping section, the cavity quickly evolves towards the corresponding point with a similar value of $\tau_b/N$ along the upsloping part.

We have compared our method with those from \citet{gagliardini2007} and \citet{stubblefield2021}. When computing the points along the steady sliding law in Figure \ref{fig:sliding} with our method, the number of time steps required to converge to a steady state can become very large (of order 1000) in the highly cavitated stages along the downsloping region of the curve. This is due to very small scale oscillations that travel across the cavity but have a significant effect on the calculated values of $\tau_b$ due to the stress singularity at the reattachment point. Contrastingly, when using the method from \citet{gagliardini2007}, these oscillations seem to dampen and the method can converge in about 100 iterations for highly cavitated states. We speculate that this is due to the use of a numerical stabilisation in Elmer when solving the advection equation, see (11) in \citet{gagliardini2008}. Despite this difference in computational times, we find that the basal stress computations carried out with our method appear to be more accurate due to the exact enforcement of discrete contact conditions (see Figure \ref{fig:steady-state} and compare with \citet[Figure 1]{gagliardini2007}). 

We also attempted to solve the subglacial cavity problem with the numerical method from \citet{stubblefield2021}. Although this method is suitable for solving the grounding line and subglacial lake problems presented in that reference, we found that it was not able to evolve the cavity correctly. This method solves the variational inequality via a penalty method, in which the constraint $\bu\cdot\bn \leq 0$ in \eqref{eq:convex_set} is enforced pointwise at the discrete level on the attached region via a penalisation term. However, due to the discretisation of the advection equation used in \citet{stubblefield2021}, enforcing this condition does not guarantee that the term $\bu \cdot \left(\partial\theta/\partial x, -1 \right)$ is exactly zero where $\bu\cdot\bn = 0$ at the discrete level. Moreover, according to our computations, when applying this penalty method to the subglacial cavity problem, the computed values of $\bu_h\cdot\bn$ along the attached region appear to be much less accurate than when using the method presented here, as shown in Appendix \ref{app:1}. We also remark that the method in \citet{stubblefield2021} does not include any stabilisation of the advection equation. As demonstrated in Appendix \ref{app:2}, without stabilisation numerical oscillations appear along the cavity roof after some time steps. For these reasons, we expect that the numerical method from \citet{stubblefield2021} should be modified for the subglacial cavity problem.

In the numerical method introduced in this paper, we are essentially solving a Stokes variational inequality by enforcing the average values of $\bu_h\cdot\bn$ to be less than or equal to 0 along the attached region. When formulated as a mixed problem, we arrive at system \eqref{eq:mixed}. However, many other possibilities exist at the discrete level; for example, we could enforce that the projection of $\bu_h\cdot \bn$ onto a space of piecewise linear functions along the attached region of the boundary be less than or equal to zero. Alternatively, we could also enforce this constraint to hold for the values of $\bu_h\cdot \bn$ at the midpoints of the edges. We are currently exploring these different discretisations in combination with different time stepping schemes for the advection equation \eqref{eq:advection}. We believe that suitable combinations of these methods could yield further improvements in terms of the speed and robustness of the algorithm. We are also exploring extensions of this method to three-dimensional problems; in particular, generalising the scheme for the advection equation is not obvious and requires further consideration.

\backsection[Acknowledgements]{We acknowledge helpful discussions with O.~Gagliardini and A.~Stubblefield.}

\backsection[Funding]{
This work was supported by the Engineering and Physical Sciences Research Council (PEF, grant numbers EP/V001493/1 and EP/R029423/1). GGdD was supported by the University of Oxford Mathematical Institute Graduate Scholarship.
}

\backsection[Declaration of interests]{The authors report no conflict of interest.}

\backsection[Data availability statement]{The code for the viscous contact solver presented in Algorithm \ref{alg:1} is openly available in the archived repository \citet{zenodo2} and in \url{https://bitbucket.org/gonzalogddiego/subglacialcavitysolver_2021/src/master/}. The solver is implemented in Firedrake \citet{rathgeber2016}, using the version available at \citet{zenodo1}. }

\backsection[Author ORCID]{G.~G.~de Diego https://orcid.org/0000-0002-9896-024X; P.~E.~Farrell, https://orcid.org/0000-0002-1241-7060; I.~J.~Hewitt https://orcid.org/0000-0002-9167-6481}

\appendix

\section{Additional numerical tests}

\subsection{Comparison with a penalty method}\label{app:1}

The method introduced in \citet{stubblefield2021} solves the variational inequality \eqref{eq:vi} by penalising the values of $\bu\cdot\bn$ which are positive on $\Gamma_a$. When enforcing the Neumann boundary condition \eqref{eq:bc-Neumann}, this is achieved by solving the variational equality 
\begin{subequations}\label{eq:stokes-penalty}
\begin{align}
	a(\bu,\bv) - b(p, \bv) + \frac{1}{\delta} \int_{\Gamma_a} \left(\bu\cdot\bn + |\bu\cdot\bn|\right) \left(\bv\cdot\bn\right)\,\dd s &= f(\bv) && \forall \bv\in V,\label{eq:stokes-penalty-1}\\
	b(q,\bu) &= 0 && \forall q\in Q,
\end{align}
\end{subequations}
where $\delta > 0$ is the penalty parameter. The finite element method in \citet{stubblefield2021} solves \eqref{eq:stokes-penalty} with Taylor-Hood elements. As $\delta \to 0$, the penalty term $\int_{\Gamma_a} \left(\bu\cdot\bn + |\bu\cdot\bn|\right) \left(\bu\cdot\bn\right)\,\dd s$ also tends to 0, which means that the contact condition $\bu\cdot \bn\leq 0$ is enforced on $\Gamma_a$. We remark that in the limit $\delta \to 0$, this method is equivalent to a Lagrange multiplier method in which the discrete contact boundary conditions set $(\bu_h\cdot \bn)(\bx)\leq 0$ for all $\bx\in\Gamma_a$, as opposed to forcing the average values of $\bu_h\cdot \bn$ to be non-positive as we do in our method.

Additionally, in \citet{stubblefield2021}, the advection equation \eqref{eq:advection-1} is solved explicitly with
\begin{equation}\label{eq:adv-stubblefield}
	\frac{\theta_i^{k+1} - \theta_i^k}{\Delta t} = - u_i \left(\frac{\partial \theta_h^k}{\partial x}\right)_i + v_i,
\end{equation}
where $(u_i,v_i)$ are the values of $\bu_h$ at the vertices in $\Gamma_a$ and $\partial \theta_h^k/\partial x$ the projection of the (distributional) derivative of $\theta_h$ onto the space of piecewise linear functions on $\Gamma_a$.

\begin{figure}
	\centering
	\includegraphics[width=\textwidth]{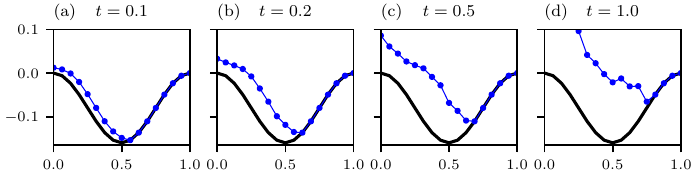}
	\includegraphics[width=\textwidth]{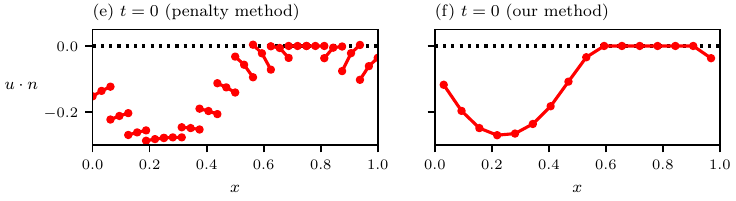}
	\caption{(a)-(d) Evolution of a subglacial cavity with the method from \citet{stubblefield2021}. (e) Values of $\bu_h\cdot\bn$ at the endpoints and midpoint of each edge computed with the penalty method from \citet{stubblefield2021} after a single timestep. (f) Average values of $\bu_h\cdot\bn$ at each edge computed after a single timestep with the method presented in this work. We set $r = 0.08$, $N = 1$, $u_i = 1$, $\Delta t = 0.01$, $n = 1$ and $A = 0.5$.}
	\label{fig:app-1}
\end{figure}

In Section \ref{subsubsec:advection}, we explain the need for a compatibility between the discrete contact boundary conditions and the discretisation of the advection equation \eqref{eq:advection-1} in order to have a robust numerical method for viscous contact problems. One might expect that solving \eqref{eq:stokes-penalty} directly should make the design of a compatible discretisation of the advection equation simpler, because we are enforcing $\bu\cdot \bn\leq 0$ on $\Gamma_a$ in a pointwise manner at the discrete level. However, our numerical tests on the subglacial cavity problem indicate that this is not the case. This appears to be because the discretisation of the advection equation in \eqref{eq:adv-stubblefield} is not compatible with the discrete contact boundary conditions in the sense that $\bu_h\cdot \bn(\bx) = 0$ for all $\bx\in\Gamma_a$ does not necessary imply that 
\begin{equation}
	- u_i \left(\frac{\partial \theta_h^k}{\partial x}\right)_i + v_i = 0
\end{equation}
on a vertex $(x_i,\theta_i)\in\Gamma_a$. Moreover, the values of $\bu_h\cdot\bn$ computed with the penalty method appear to be less accurate then the corresponding values computed with the method presented in this work. Due to this lack of accuracy, we no longer have a clear distinction between regions that detach ( $\bu\cdot\bn < 0$) and those that stay attached ($\bu\cdot \bn = 0$). 

This latter issue is illustrated in Figure \ref{fig:app-1}. In this test, we evolve the cavity by solving the variational inequality equipped with the Dirichlet boundary conditions \eqref{eq:bc-Dirichlet} with the penalty method described above with $r = 0.08$ over a mesh with 16 elements along the lower boundary. We set $N = 1$, $u_i = 1$, $\Delta t = 0.01$, $n = 1$, $A = 0.5$ and $\delta = 10^{-6}$ in \eqref{eq:stokes-eq-1}. Additionally, in order to avoid detachment, we increase the tolerance in \eqref{eq:criterion_contact} to $\mathrm{tol} = 10^{-3}$. In Figures \ref{fig:app-1}(a)-(d) we present the evolution of the cavity from a fully attached state when computed our implementation of the method from \citet{stubblefield2021}. Our numerical results show that by $t = 0.5$ the cavity shape is deforming excessively. In Figures \ref{fig:app-1}(e) and (f) we also compare the computed values of $\bu_h\cdot\bn$ in the first time step with the penalty method from \citet{stubblefield2021} and with our method, respectively. As seen in Figure \ref{fig:app-1}(e), we have points where $\bu_h\cdot\bn > 0$ in between points where $\bu_h\cdot\bn = 0$.

We have carried out further computations with more refined meshes, and these indicate that the lack of accuracy in $\bu_h\cdot\bn$ persists when using the penalty method. We would like to emphasise that this numerical issue does not arise when solving the problems presented in \citet{stubblefield2021}. It must therefore be due to some particularity of the subglacial cavity problem.

\subsection{On the use of upwinding when solving the advection equation}\label{app:2}

In Section \ref{subsubsec:advection} we state that the advection equation \eqref{eq:advection-1} should be stabilised with e.g.~an upwinding method. Since no stabilisation is required for the numerical tests solved in \citet{stubblefield2021}, here we demonstrate that without upwinding, the subglacial cavity problem quickly becomes unstable. 

\begin{figure}
	\centering
	\includegraphics[width=\textwidth]{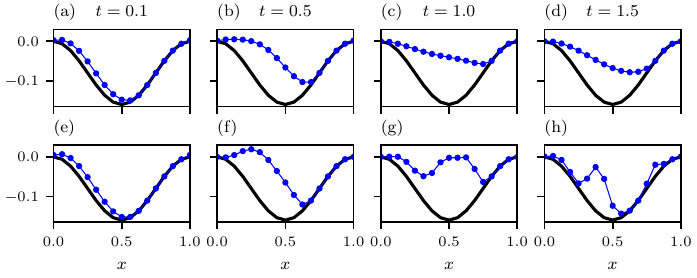}
	\caption{Evolution of the cavity roof with upwinding ((a)-(d)) and without upwinding ((e)-(h)). For this numerical test, we set $r = 0.08$, $N = 1$, $u_i = 1$, $\Delta t = 0.05$, $n = 1$ and $A = 0.5$.}
	\label{fig:app-2}
\end{figure}

In order to show this, we solve the same problem as in Appendix \ref{app:1} with a new method that is identical to the one presented in this paper except for the fact that the advection equation is discretised as
\begin{equation}\label{eq:discrete-advection-noupwind}
	\frac{\theta_i^{k+1} - \theta_i^k}{\Delta t} = \left( \left(\dfrac{\theta^k_i - \theta^k_{i-1}}{x_i - x_{i-1}}\right)^2 + 1 \right)^{1/2} \left( \frac{u^k_{n,i} + u^k_{n,i+1}}{2} \right).
\end{equation}
instead of \eqref{eq:discrete-advection}. That is, instead of taking the average values of $\bu_h\cdot\bn$ along the edge immediately upstream each node, the average of the edges upstream and downstream is taken. In this way, we still have compatibility with the discrete boundary conditions but no upwinding. In Figure \ref{fig:app-2} we can see the evolution of the cavity with and without upwinding. Oscillations clearly arise in the cavity roof at $t = 1$ when no upwinding is used.

\section{Linearised solution}\label{app:linear}

For a Newtonian flow ($n = 1$) and small amplitude topography, the theory of \citet{fowler1986} and \citet{schoof2005} can be used to find solutions with which the numerical calculations can be compared.  The method involves linearisation of the boundary conditions and use of complex variables to solve a Riemann-Hilbert problem for the velocities and stresses.  We summarise the result for the particular case when $b = r L \cos\left(2\pi x/L\right)$.  

It is convenient to parameterise the solution in terms of the scaled cavity end points $c$ and $d$, such that the cavity occupies the region  $d<\hat{x}<c+1$, where $\hat{x} = x/L$.  The velocity is $\mathbf{u} = (u,v)$ with $u \approx u_b$, and on the cavitated region the vertical velocity $v$ satisfies
\begin{equation}
	\begin{split}
		\frac{\partial v}{\partial x} = & -\frac{4\pi^2 r u_b}{L} \cos\left( \frac{2\pi x}{L}\right) \\ 
		& + \frac{4\pi^2 r u_b}{L} \left| \frac{\sin\pi(d-\hat{x})}{\sin\pi(\hat{x}-c)}\right|^{1/2}\left[\cos \pi \left(2\hat{x}+\half(d-c)\right) -\sin\pi(c+d)\sin\frac{\pi}{2}(d-c) \right].
	\end{split}
\end{equation}
The linearised steady kinematic condition for the cavity roof is $u_b \frac{\partial \theta}{\partial x} = v$, and integrating this subject to the conditions that $\theta = b$ at the cavity end points provides a constraint between $c$ and $d$.  In addition, $c$ and $d$ are related to the effective pressure by
\begin{equation}
	N = \frac{8 \pi^2 r \eta u_b}{L} \cos\frac{\pi}{2}(3d+c) \sin \frac{\pi}{2}(d-c).
\end{equation}
Thus, for given values of $N$ and $u_b$, these two constraints determine the end points $c$ and $d$.  Further, the normal stress on the contact region $c < \hat{x} < d$ is given by
\begin{equation}
	\begin{split}
	& \sigma_{nn} +p_w = \\
	&  - \frac{8\pi^2 r \eta u_b}{L} \left| \frac{\sin\pi(d-\hat{x})}{\sin\pi(\hat{x}-c)}\right|^{1/2}\left[\cos\pi \left(2\hat{x}+\half(d-c)\right) -\sin\pi(c+d)\sin\frac{\pi}{2}(d-c) \right],
	\end{split}
\end{equation}
and the integral in \eqref{eq:basal} then gives the basal shear stress as
\begin{equation}
	\begin{split}
	\tau_b = \frac{\pi^3 r^2 \eta u_b}{L} \left[  5 - \sin 2\pi(c+d) \sin 2\pi(d-c)- \cos2\pi(d-c)-4\cos\pi(d-c)  \right. \\ \left. -\cos\pi(3c+d)+\cos\pi(c+3d) \right].
	\end{split}
\end{equation}
This cavitated solution requires $N < 8 \pi^2 r \eta u_b/L$.  Otherwise, there is no cavity and we have $\tau_b = 8\pi^3 r^2 \eta u_b/L$.

\bibliographystyle{jfm}
\bibliography{bibliography}

\end{document}